\documentclass{article}
\usepackage{amsmath}
\numberwithin{equation}{section}
\usepackage{indentfirst}

\begin{document}

\newtheorem{thm}{Theorem}[section]
\newtheorem{cor}[thm]{Corollary}
\newtheorem{prop}[thm]{Proposition}
\newtheorem{conj}[thm]{Conjecture}
\newtheorem{lem}[thm]{Lemma}
\newtheorem{Def}[thm]{Definition}
\newtheorem{rem}[thm]{Remark}
\newtheorem{prob}[thm]{Problem}
\newtheorem{ex}{Example}[section]

\newcommand{\be}{\begin{equation}}
\newcommand{\ee}{\end{equation}}
\newcommand{\ben}{\begin{enumerate}}
\newcommand{\een}{\end{enumerate}}
\newcommand{\beq}{\begin{eqnarray}}
\newcommand{\eeq}{\end{eqnarray}}
\newcommand{\beqn}{\begin{eqnarray*}}
\newcommand{\eeqn}{\end{eqnarray*}}
\newcommand{\bei}{\begin{itemize}}
\newcommand{\eei}{\end{itemize}}

\newcommand{\pa}{{\partial}}
\newcommand{\V}{{\rm V}}
\newcommand{\R}{{\rm R}}
\newcommand{\e}{{\epsilon}}
\newcommand{\tomega}{\tilde{\omega}}
\newcommand{\tOmega}{\tilde{Omega}}
\newcommand{\tR}{\tilde{R}}
\newcommand{\tB}{\tilde{B}}
\newcommand{\tGamma}{\tilde{\Gamma}}
\newcommand{\fa}{f_{\alpha}}
\newcommand{\fb}{f_{\beta}}
\newcommand{\faa}{f_{\alpha\alpha}}
\newcommand{\faaa}{f_{\alpha\alpha\alpha}}
\newcommand{\fab}{f_{\alpha\beta}}
\newcommand{\fabb}{f_{\alpha\beta\beta}}
\newcommand{\fbb}{f_{\beta\beta}}
\newcommand{\fbbb}{f_{\beta\beta\beta}}
\newcommand{\faab}{f_{\alpha\alpha\beta}}

\newcommand{\pxi}{ {\pa \over \pa x^i}}
\newcommand{\pxj}{ {\pa \over \pa x^j}}
\newcommand{\pxk}{ {\pa \over \pa x^k}}
\newcommand{\pyi}{ {\pa \over \pa y^i}}
\newcommand{\pyj}{ {\pa \over \pa y^j}}
\newcommand{\pyk}{ {\pa \over \pa y^k}}
\newcommand{\dxi}{{\delta \over \delta x^i}}
\newcommand{\dxj}{{\delta \over \delta x^j}}
\newcommand{\dxk}{{\delta \over \delta x^k}}

\newcommand{\px}{{\pa \over \pa x}}
\newcommand{\py}{{\pa \over \pa y}}
\newcommand{\pt}{{\pa \over \pa t}}
\newcommand{\ps}{{\pa \over \pa s}}
\newcommand{\pvi}{{\pa \over \pa v^i}}
\newcommand{\ty}{\tilde{y}}
\newcommand{\bGamma}{\bar{\Gamma}}

\font\BBb=msbm10 at 11pt
\newcommand{\Bbb}[1]{\mbox{\BBb #1}}

\newcommand{\qed}{\hspace*{\fill}Q.E.D.}  

\title{Some fundamental problems in global Finsler geometry}
\author{ Xinyue Cheng\footnote{supported by the National Natural Science Foundation of China (11871126) and the Science Foundation of Chongqing Normal University (17XLB022)}}


\maketitle

\begin{abstract}
The geometry and analysis on Finsler manifolds is a very important part of Finsler geometry. In this article, we introduce some important and fundamental topics in global Finsler geometry and discuss the related properties and the relationships in them. In particular, we optimize and improve the various definitions of Lie derivatives on Finsler manifolds. We also characterize the gradient vector fields  and obtain a gradient estimate  for any smooth function on a Randers manifold.\\
{\bf Keywords:} dual Finsler metric; gradient vector field; Finsler Laplacian; eigenvalue; Hessian; Lie derivative; gradient estimate. \\
{\bf  MR(2000) Subject Classification:} 53B40,  53C60
\end{abstract}

\section{Preliminaries}

Let $M$ be a connected manifold of dimension $n$ and $\pi : TM_{0} \rightarrow M$ be the natural projective map, where $TM_{0}:=TM\setminus \{0\}$. $\pi$ pulls back $TM$ to a vector bundle $\pi ^{*}TM$ over $TM_{0}$. The fiber at a point $(x,y)\in TM_{0}$ is defined by
\[
\pi ^{*}TM|_{(x,y)}:=\left\{(x,y, v) \ | \ v\in T_{x}M\right\}\cong T_{x}M.
\]
In other words, $\pi ^{*}TM|_{(x,y)}$ is just a copy of $T_{x}M$. Similarly, we define the pull-back cotangent bundle $\pi ^{*}T^{*}M$ whose fiber at $(x,y)$ is a copy of $T^{*}_{x}M$. That is,
\[
\pi ^{*}T^{*}M|_{(x,y)}:=\left\{(x,y, \theta ) \ | \ \theta\in T_{x}^{*}M\right\}\cong T^{*}_{x}M.
\]
$\pi ^{*}T^{*}M$ can be viewed as the dual vector bundle of $\pi ^{*}TM$ by setting
\[
(x, y, \theta)(x, y, v):= \theta (v), \ \ \ \theta \in T^{*}_{x}M, \ v \in T_{x}M.
\]

Take a standard local coordinate system $(x^{i}, y^{i})$ in $TM$. Let $\{\frac{\pa}{\pa x^{i}}, \frac{\pa}{\pa y^{i}}\}$ and $\{dx^{i}, dy^{i}\}$ be the natural local frame and coframe for $T(TM_{0})$ and  $T^{*}(TM_{0})$ respectively. Let
\[
\pa _{i}:=\left(x, y, \frac{\pa}{\pa x^{i}}|_{x}\right).
\]
Then $\{\pa _{i}\}$ is a local frame for $\pi ^{*}TM$. Dually, put
\[
dx^{i}:=\left(x, y,  dx^{i}|_{x}\right).
\]
Then $\{dx^{i}\}$ is a local coframe for  $\pi ^{*}T^{*}M$.

The vertical tangent bundle of $M$ is defined by $VTM:= span\{\frac{\pa}{\pa y^{i}}\}$. $VTM$ is a well-defined subbundle of $T(TM_{0})$ and we can obtain a decomposition $T(TM_{0})=\pi ^{*}TM \oplus VTM$.

For a Finsler manifold $(M, F)$, let
\be
G = y^{i}\frac{\pa}{\pa x^{i}}-2 G^{i}\frac{\pa}{\pa y^{i}}, \label{Spray}
\ee
where $G^{i}=G^{i}(x, y)$ are defined by
\[
G^{i}=\frac{1}{4}\left\{[F^2]_{x^{k}y^{l}}y^{k}-[F^2]_{x^{l}}\right\}.
\]
We call $G$ the spray induced by $F$ and $G^i$ the spray coefficients of $F$. Define $N^{i}_{j}:=\frac{\pa G^{i}}{\pa y^{j}}$ and let
\be
\frac{\delta}{\delta x^{i}}:= \frac{\pa}{\pa x^{i}}- N_{i}^{j}\frac{\pa}{\pa y^{j}}.
\ee
Then $\left\{\frac{\delta}{\delta x^{i}}, \frac{\pa}{\pa y^i}\right\}$ form a local frame for $T(TM)$. Further, $HTM:=span \left\{\frac{\delta}{\delta x^{i}}\right\}$ is a well-defined subbundle of $T(TM_{0})$ and is called the horizontal tangent bundle of $M$. Then we obtain a decomposition for  $T(TM_{0})$, $T(TM_{0})=HTM \oplus VTM$.

The following maps are natural and are important for our discussions below.
\ben
\item[{\rm (1)}] Define a vector bundle map $\rho : T(TM_{0})\rightarrow \pi ^{*}TM$ by
\be
\rho \left(\frac{\pa}{\pa x^{i}}|_{(x,y)}\right)=\pa _{i}, \ \ \ \rho \left(\frac{\pa}{\pa y^{i}}|_{(x,y)}\right)=0. \label{Defofrho}
\ee
It is clear that $\ker \rho = VTM$.
\item[{\rm (2)}] Define a linear map ${\cal H}: \pi ^{*}TM \rightarrow HTM$ with the following properties
\be
{\cal H}(\pa _{i}):=\frac{\delta}{\delta x^{i}}. \label{defofH}
\ee
Obviously, ${\cal H}$ is an isomorphism.
\een

\section{Gradient vector fields and Laplacian on Finsler manifolds}

Let $M$ be an $n$-dimensional manifold. A Finsler metric $F$ on $M$ is a non-negative function on $TM$ such that $F$ is $C^{\infty}$ on $TM\backslash \{0\}$ and the restriction $F_{x}:=F|_{T_{x}M}$ is a Minkowski function on $T_{x}M$ for all $x\in M$. For Finsler metric $F$ on $M$, there is a Finsler co-metric $F^{*}$ on $M$ which is non-negative function on the cotangent bundle $T^{*}M$ given by
\be
F^{*}(x, \xi):=\sup\limits_{y\in T_{x}M\setminus \{0\}} \frac{\xi (y)}{F(x,y)}, \ \ \forall \xi \in T^{*}_{x}M.
\ee
We call $F^{*}$ the dual Finsler metric of $F$. Finsler metric $F$ and its dual Finsler metric $F^{*}$ satisfy the following relation.

\begin{lem}{\rm (Lemma 3.1.1, \cite{shen1})}\label{shen311} Let $F$ be a Finsler metric on $M$ and $F^{*}$ its dual Finsler metric. For any vector $y\in T_{x}M\setminus \{0\}$, $x\in M$, the covector $\xi =g_{y}(y, \cdot)\in T^{*}_{x}M$ satisfies
\be
F(x,y)=F^{*}(x, \xi)=\frac{\xi (y)}{F(x,y)}. \label{shenF311}
\ee
Conversely, for any covector $\xi \in T_{x}^{*}M\setminus \{0\}$, there exists a unique vector $y\in T_{x}M\setminus \{0\}$ such that $\xi =g_{y}(y, \cdot)\in T^{*}_{x}M$.
\end{lem}

Naturally, by Lemma \ref{shen311}, we define a map ${\cal L}: TM \rightarrow T^{*}M$ by
\[
{\cal L}(y):=\left\{
\begin{array}{ll}
g_{y}(y, \cdot), & y\neq 0, \\
0, & y=0.
\end{array} \right.
\]
It follows from (\ref{shenF311}) that
\be
F(x,y)=F^{*}(x, {\cal L}(y)).
\ee
Thus ${\cal L}$ is a norm-preserving transformation. We call ${\cal L}$ the Legendre transformation on Finsler manifold $(M, F)$.

Take a basis $\{{\bf b}_{i}\}^{n}_{i=1}$ for $TM$ and its dual basis $\{\theta ^{i}\}_{i=1}^{n}$ for $T^{*}M$. Express
\[
\xi ={\cal L}(y)=\xi _{i}\theta ^{i}=g_{ij}(x,y)y^{j}\theta ^{i},
\]
where $g_{ij}(x,y):=\frac{1}{2}\left[F^2\right]_{y^{i}y^{j}}(x,y)$. The Jacobian of ${\cal L}$ is given by
\[
\frac{\pa \xi _{i}}{\pa y^{j}}=g_{ij}(x,y).
\]
Thus ${\cal L}$ is a diffeomorphism from $TM\setminus\{0\}$ onto $T^{*}M\setminus\{0\}$.  Let
\be
g^{*kl}(x,\xi):=\frac{1}{2}\left[F^{*2}\right]_{\xi _{k}\xi_{l}}(x,\xi).
\ee
For any $\xi ={\cal L}(y)$, differentiating $F^{2}(x,y)=F^{*2}(x,{\cal L}(y))$ with respect to $y^i$ yields
\be
\frac{1}{2}\left[F^2\right]_{y^{i}}(x,y)=\frac{1}{2}\left[F^{*2}\right]_{\xi _{k}}(x,\xi)g_{ik}(x,y), \label{shen3121}
\ee
which implies
\be
g^{*kl}(x,\xi)\xi _{l}=\frac{1}{2}\left[F^{*2}\right]_{\xi _{k}}(x,\xi)=\frac{1}{2}g^{ik}(x,y)\left[F^2\right]_{y^{i}}(x,y)=y^{k}. \label{dualY}
\ee
Then, it is clear that
\[
g^{*kl}\xi _{l}\frac{\pa g_{ik}}{\pa y^{j}}=y^{k}\frac{\pa g_{ik}}{\pa y^{j}}=0.
\]
Differentiating (\ref{shen3121}) with respect to $y^{j}$ gives
\beqn
g_{ij}(x,y)&=& g^{*kl}(x,\xi)g_{ik}(x,y)g_{jl}(x,y)+g^{*kl}(x,\xi)\xi _{l}\frac{\pa g_{ik}}{\pa y^{j}}(x,y) \\
&=& g^{*kl}(x,\xi)g_{ik}(x,y)g_{jl}(x,y).
\eeqn
Therefore, we get
\be
g^{*kl}(x,\xi)=g^{kl}(x,y).  \label{Fdual}
\ee

\vskip 3mm

Given a smooth function $f$ on $M$, the differential $df_{x}$ at any point $x\in M$,
\[
df_{x}=\frac{\pa f}{\pa x^i}(x)dx^{i}
\]
is a linear function on $T_{x}M$. We define the gradient vector $\nabla f(x)$ of $f$ at $x\in M$ by $\nabla f(x):={\cal L}^{-1}\left(df(x)\right)\in T_{x}M$. In a local coordinate system, by (\ref{dualY}) and (\ref{Fdual}), we can express $\nabla f$ as
\be
\nabla f(x)=\left\{
\begin{array}{ll}
g^{ij}(x,\nabla f)\frac{\pa f}{\pa x^i}\frac{\pa}{\pa x^j}, & x\in M_{f},  \label{gradientV}\\
0, & x \in M\setminus M_{f},
\end{array}
\right.
\ee
where $M_{f}=\{x\in M| df (x)\neq 0\}.$ Further, by Lemma \ref{shen311}, we have the following
\be
df_{x}(v)=g_{\nabla f_{x}}(\nabla f_{x}, v), \ \ \ \forall v\in T_{x}M
\ee
and
\be
F(x, \nabla f_{x})=F^{*}(x, df_{x})=\frac{df_{x} (\nabla f_{x})}{F(x,\nabla f_{x})}.
\ee

\vskip 4mm

By definition, a smooth measure $\mu$ on $M$ is a measure locally given by a smooth $n$-form
\[
d\mu =\sigma (x)dx^{1}\cdots dx^{n}.
\]
The restriction $\mu _{x}$ of $\mu$ to $T_{x}M$ is a Haar measure on $T_{x}M$. For every Finsler manifold $(M, F)$, there are several associated measures, including Busemann-Hausdorff measure $\mu _{BH}$ and Holmes-Thompson measure $\mu _{HT}$. A Finsler manifold $(M,F)$ equipped with a smooth measure $\mu$ is called a Finsler measure space and denoted by $(M, F, d\mu)$.

Let us consider an oriented manifold $M$ equipped with a measure $\mu$. We can view $d\mu$ as an $n$-form (volume form) on $M$. Let $X$ be a vector field on $M$. Define an $(n-1)$-form $X\rfloor d\mu$ on $M$ by
\[
X\rfloor d\mu (X_{2}, \cdots , X_{n}):=d\mu (X, X_{2}, \cdots , X_{n}).
\]
Define
\be
d(X\rfloor d\mu)=div (X)d\mu. \label{divergence}
\ee
We call $div(X)$ the divergence of $X$. Clearly, $div (X)$ depends only on the volume form $d\mu$. In a local coordinate system $(x^i)$, express $d\mu =\sigma (x)dx^{1}\cdots dx^{n}$. Then for a vector field $X=X^{i}\frac{\pa}{\pa x^i}$ on $M$,
\be
div(X)=\frac{1}{\sigma}\frac{\pa}{\pa x^i}\left(\sigma X^{i}\right)=\frac{\pa X^i}{\pa x^i}+\frac{X^i}{\sigma}\frac{\pa \sigma}{\pa x^i}.
\ee
Applying the Stokes theorem to $\eta =X\rfloor d\mu$, we obtain
\beq
\int _{M} div(X)d\mu =\int _{M} d(X\rfloor d\mu)=0, & if \ \pa M=\emptyset , \label{divergence1}\\
\int _{M} div(X)d\mu =\int _{M} d(X\rfloor d\mu)=\int _{\pa M}X\rfloor d\mu , & if \ \pa M \neq \emptyset . \label{divergence2}
\eeq

\vskip 3mm

Now we introduce the Laplacian on a Finsler  measure space $(M, F, d\mu)$. There are some different definitions on Laplacian in Finsler geometry (e.g. see \cite{AL}\cite{HeSh}\cite{BTh}). The following definition is from \cite{GeShen}.

Given a smooth measure $d\mu =\sigma (x)dx^{1}\cdots dx^{n}$ and $C^k$ ($k\geq 2$) function $f$ on $M$, define the Finsler Laplacian $\triangle f$ of $f$ by
\be
\triangle f:=div(\nabla f). \label{laplace}
\ee
By (\ref{gradientV}), the Laplacian of $f$ is expressed by
\be
\triangle f= \frac{1}{\sigma}\frac{\pa}{\pa x^i}\left(\sigma\nabla ^{i}f\right)=\frac{1}{\sigma}\frac{\pa}{\pa x^i}\left(\sigma g^{ij}(x,\nabla f)\frac{\pa f}{\pa x^j}\right), \label{laplace2}
\ee
where $\nabla ^{i}f:= g^{ij}(x,\nabla f)\frac{\pa f}{\pa x^j}= g^{*ij}(x, df)\frac{\pa f}{\pa x^j}$. From (\ref{laplace2}), Finsler Laplacian is a nonlinear elliptic differential operator of the second order.

\begin{rem} {\rm The following are some remarks on Finsler Laplacian.
\ben
\item[{\rm (i)}] As we know, for any smooth function $\varphi$ on $M$,
\[
div(\varphi \nabla f)=\varphi \triangle f+d\varphi (\nabla f).
\]
If $\pa M =\emptyset$, applying the divergence formula (\ref{divergence1}) to the above identity yields
\be
\int_{M}\varphi \triangle fd\mu =-\int _{M}d\varphi (\nabla f)d\mu. \label{shen145}
\ee
Actually, (\ref{shen145}) gives the definition of non-linear Laplacian $\triangle f$ on the whole $M$ in the distribution sense.
\item[{\rm (ii)}] Let $\mu =e^{\rho}\mu _{F}$ and $\triangle _{F}$ denote the Laplacian associated with the Finsler measure $\mu _{F}$. Then
\[
\triangle f=\triangle _{F}f+d\rho (\nabla f).
\]
Note that $\triangle f$ and $\triangle _{F}f$ are the divergences of the gradient $\nabla f$ with respect to $\mu$ and $\mu _{F}$ respectively.
\item[{\rm (iii)}] The Finsler $p$-Laplacian $\triangle _{p}f$ of $f$ is formally defined by
\[
\triangle _{p}f:=div\left[F^{p-2}(x, \nabla f)\nabla f\right].
\]
In the distribution sense, the definition of Finsler $p$-Laplacian $\triangle _{p}f$ is given by the following identity
\[
\int _{M}\varphi \triangle _{p}f d\mu =-\int _{M}F^{p-2}(x,\nabla f)d\varphi (\nabla f)d\mu, \ \ \forall \varphi \in C^{\infty}_{0}(M).
\]
When $p=2$, $\triangle _{p}$ is exactly the usual Finsler Laplacian.
\een
}
\end{rem}

\section{Energy functionals and eigenvalues}

The variational problem of the canonical energy functional also gives rise to the Laplacian. Let $H^{1}$ denote the Hilbert space of all $L^2$ functions  $f\in C^{\infty}(M)$ such that $df \in L^2$. Denote by $H^{1}_{0}$ the space of functions $u\in H^{1}$ with $\int _{M}ud\mu =0$ if $\pa M=\emptyset $ and with $u_{|\pa M}=0$ if $\pa M \neq \emptyset $. The canonical energy functional ${\cal E}$ on $H^{1}_{0}$ is defined by
\[
{\cal E}(u):=\frac{\int _{M}\left[F^{*}(x,du)\right]^{2}d\mu}{\int _{M}u^{2}d\mu}.
\]
For functions $u, \varphi \in H^{1}_{0}$, by (\ref{dualY}), we have
\[
\frac{d}{d\varepsilon}\left[F^{*2}(x,du +\varepsilon d\varphi)\right]|_{\varepsilon =0}=\frac{\pa [F^{*2}]}{\pa \xi _{i}}(x, du)\frac{\pa \varphi}{\pa x^{i}}=2\nabla ^{i}u(x, du)\frac{\pa \varphi }{\pa x^{i}}=2d\varphi (\nabla u).
\]
Thus, for any $u\in H^{1}_{0}$ with $\int _{M}u^{2}d\mu=1$,
\be
d_{u}{\cal E}(\varphi)=2\int d\varphi(\nabla u)d\mu-2\lambda \int u\varphi \ d\mu , \ \ \ \forall \varphi \in H^{1}_{0}, \label{shen148}
\ee
where $\lambda ={\cal E} (u)> 0$. From (\ref{shen145}), we can rewrite (\ref{shen148}) as follows
\be
\frac{1}{2}d_{u}{\cal E}(\varphi)=-\int \left[\triangle u+\lambda u\right]\varphi \ d\mu, \ \ \ \forall \varphi \in H^{1}_{0}.
\ee
Hence, it follows that a function $u\in H^{1}_{0}$ satisfies $d_{u}{\cal E}=0$ with $\lambda ={\cal E}(u)$ if and only if
\be
\triangle u+\lambda u=0.
\ee
In  this case, $\lambda$ and $u$ are called an eigenvalue and an eigenfunction of $(M,F, d\mu)$, respectively. Thus an eigenfunction $u$ corresponding to an eigenvalue $\lambda$ satisfies the following equation
\[
\frac{1}{\sigma (x)}\frac{\pa}{\pa x^i}\left(\sigma (x)\nabla ^{i}u(x)\right)+\lambda u=0,
\]
where $\nabla ^{i}u(x)=g^{ij}(x,\nabla u)\frac{\pa u}{\pa x^j}= g^{*ij}(x, du)\frac{\pa u}{\pa x^j}.$

Denote by ${\cal E}_{\lambda}$ the union of the zero function and the set of all eigenfunctions corresponding to $\lambda$. We call ${\cal E}_{\lambda}$ the eigencone corresponding to $\lambda$.

Assume that $M$ is compact without boundary. Let
\be
\lambda _{1}(M):=\inf \limits_{u\in C^{\infty}(M)}\frac{\int _{M}\left[F^{*}(x, du)\right]^{2}d\mu}{\inf _{\lambda \in R}\int _{M}|u-\lambda|^{2}d\mu}. \label{firste}
\ee
In \cite{GeShen} and \cite{shen1}, we can find the proof on the fact that $\lambda _{1}(M)$ is the minimum of the energy functional ${\cal E}$. In the following, we give a different proof for this fact.
\begin{prop} \  $\lambda _{1}:=\lambda _{1}(M)$ is the smallest eigenvalue of $(M,F, d\mu)$, that is,  $\lambda _{1}=\inf _{u\in H^{1}_{0}}{\cal E}(u)$.
\end{prop}
{\bf Proof.}  Write
\beqn
\int _{M}|u-\lambda |^2 d\mu &=& \int _{M}u^{2}d\mu -2\lambda \int _{M}ud\mu +\lambda ^{2}\int _{M}d\mu \\
     &:=& a-2\lambda b+\lambda ^{2}c,
\eeqn
where $a:=\int _{M}u^{2}d\mu , \ b:=\int _{M}ud\mu , c:=\int _{M}d\mu =\mu (M)$.  Let $f(\lambda):= a-2\lambda b+\lambda ^{2}c$. By $f'(\lambda)=-2b+2c\lambda$, we have the following
\beqn
&& \inf\limits_{\lambda\in R}\int _{M}|u-\lambda |^2 d\mu= \left.\left(\int _{M}|u-\lambda |^2 d\mu \right)\right|_{\lambda =\frac{b}{c}} \\
&&=a-2\mu(M)^{-1}b^{2}+\mu (M)^{-1}b^{2} \\
&&=\int _{M}u^{2}d\mu -\mu (M)^{-1}\left(\int _{M}ud\mu \right)^{2}\leq \int _{M}u^{2}d\mu .
\eeqn
Thus
\[
\lambda _{1}(M)=\inf \limits_{u\in H^{1}_{0}(M)}\frac{\int _{M}\left[F^{*}(x, du)\right]^{2}d\mu}{\int _{M}u^{2}d\mu}=\inf\limits_{u\in H^{1}_{0}(M)}{\cal E}(u).
\]  \qed
\vskip 2mm

We call $\lambda _{1}$ the first eigenvalue of $(M,F, d\mu)$. They are natural problems to determine the lower bound of the first (nonzero) eigenvalue of Laplacian on Finsler manifolds and to study the structure of the first eigencone for a general Finsler metric (\cite{shenshen}\cite{WX}\cite{Xia3}).

\section{Hessian }

Let $(M,F)$ be a Finsler manifold of dimension $n$ and $\pi : TM\setminus \{0\} \rightarrow M$ be the projective map. The pull-back $\pi ^{*}TM$ admits a unique linear connection, which is called the Chern connection. The Chern connection $D$ is determined by the following equations
\beq
&& D^{V}_{X}Y-D^{V}_{Y}X=[X,Y], \label{chern1}\\
&& Zg_{V}(X,Y)=g_{V}(D^{V}_{Z}X,Y)+g_{V}(X,D^{V}_{Z}Y)+C_{V}(D^{V}_{Z}V,X,Y) \label{chern2}
\eeq
for $V\in TM\setminus \{0\}$  and $X, Y, Z \in TM$, where
\[
C_{V}(X,Y,Z):=C_{ijk}(x,V)X^{i}Y^{j}Z^{k}=\frac{1}{4}\frac{\pa ^{3}F^{2}(x,V)}{\pa V^{i}\pa V^{j}\pa V^{k}}X^{i}Y^{j}Z^{k}
\]
is the Cartan tensor of $F$ and $D^{V}_{X}Y$ is the covariant derivative with respect to the reference vector $V$.

Let $(M,F)$ be a Finsler manifold. There are two ways to define the Hessian of a $C^2$ function on $M$. Let $f$ be a $C^2$ function on $M$. Firstly, the Hessian of $f$ can be defined as a map $D^{2}f:  TM\rightarrow R$ by
\be
D^{2}f(y):=\frac{d^2}{ds^2}\left(f\circ c\right)|_{s=0}, \ \ \ y\in T_{x}M,\label{hessian1}
\ee
where $c: (-\varepsilon , \varepsilon) \rightarrow M$ is the geodesic with $c(0)=x, \ \dot{c}(0)=y\in T_{x}M$ (see (\cite{shen1})).  In local coordinates,
\beq
D^{2}f(y)&=& \frac{\pa ^{2}f}{\pa x^{i}\pa x^{j}}(x)\dot{c}^{i}(0)\dot{c}^{j}(0)+\frac{\pa f}{\pa x^{i}}(x)\ddot{c}^{i}(0) \nonumber\\
&=& \frac{\pa ^{2}f}{\pa x^{i}\pa x^{j}}(x)y^{i}y^{j}-2\frac{\pa f}{\pa x^{i}}(x)G^{i}(x,y) \nonumber \\
&=& \left(\frac{\pa ^{2}f}{\pa x^{i}\pa x^{j}}(x)-\frac{\pa f}{\pa x^{m}}\Gamma ^{m}_{ij}(x,y)\right)y^{i}y^{j}.
\eeq
Here, $\Gamma ^{k}_{ij}(x,y)$ denote the Chern connection coefficients of $F$, which depends on the tangent vector $y\in T_{x}M$ usually.

There is another definition of the Hessian in Finsler geometry, by which the Hessian of a $C^2$ function $u$ on $M$ is corresponding to a symmetric matrix $\left(u_{|i|j}(x,\nabla u)\right)$, where ``$|$" denotes the horizontal covariant derivative with respect to the Chern connection of the metric. Concretely, the Hessian $\nabla ^{2}u$ of $u$ is defined by
\be
\nabla ^{2}u(X,Y):=g_{\nabla u}\left(D_{X}^{\nabla u}\nabla u, Y\right) \label{hessian2}
\ee
for any $X, Y \in TM$ (\cite{OS}\cite{WX}\cite{WuXin}). In a local coordinate system, let $X=X^{i}\frac{\pa}{\pa x^{i}}, \   Y=Y^{j}\frac{\pa}{\pa x^{j}}$. By the definition,
\[
D_{X}^{\nabla u}\nabla u =\left\{\frac{\pa (\nabla ^{i}u)}{\pa x^{j}}X^{j}+(\nabla^{k}u)\Gamma ^{i}_{jk}(x, \nabla u)X^{j}\right\}\frac{\pa}{\pa x^{i}}.
\]
Thus, we have
\beq
g_{\nabla u}\left(D_{X}^{\nabla u}\nabla u, Y\right)&=&\left(\frac{\pa (\nabla ^{i}u)}{\pa x^{j}}+(\nabla^{k}u)\Gamma ^{i}_{jk}(x, \nabla u)\right)X^{j}Y^{l}g_{il}(x,\nabla u) \nonumber\\
&=&(\nabla ^{i}u)_{|j}(x,\nabla u) g_{il}(x,\nabla u)X^{j}Y^{l}.   \label{CCD}
\eeq
Here, we have used the  facts that $\nabla ^{i}u=g^{ij}(x, \nabla u)\frac{\pa u}{\pa x^{j}}$ and
\beqn
 \frac{\pa (\nabla ^{i}u)}{\pa y^{m}}(x, \nabla u)&=&-2C_{klm}(x,\nabla u)g^{ik}(x,\nabla u)g^{jl}(x,\nabla u)\frac{\pa u}{\pa x^{j}} \\
& =& -2C_{klm}(x,\nabla u)(\nabla ^{l}u)g^{ik}(x,\nabla u)=0.
\eeqn
Further, let $(\nabla u)_{i}(x, \nabla u):=g_{ij}(x, \nabla u)(\nabla ^{j}u)$. Then
$$
(\nabla u)_{i}(x, \nabla u) =g_{ij}(x, \nabla u)g^{jk}(x,\nabla u)\frac{\pa u}{\pa x^k}=\frac{\pa u}{\pa x^{i}}=u_{|i}(x).
$$
Thus it follows (\ref{CCD}) that
\[
g_{\nabla u}\left(D_{X}^{\nabla u}\nabla u, Y\right)=u_{|i|j}(x, \nabla u)X^{i}Y^{j}.
\]
Hence, by (\ref{hessian2}), we have the following proposition.
\begin{prop}{\rm Let $u$ be  a $C^2$ function on Finsler manifold $(M, F)$. Then, for any $X, Y \in TM$,  we have
\be
\nabla ^{2}u(X,Y)=u_{|i|j}(x, \nabla u)X^{i}Y^{j}. \label{hessian3}
\ee
}
\end{prop}

It follows from (\ref{hessian3}) that the Hessian $\nabla ^{2}u$ of a $C^2$ function $u$ is determined completely by the following symmetric matrix
\be
{\rm Hess}(u): =\left(u_{|i|j}(x, \nabla u)\right).  \label{hessian}
\ee
\begin{rem}
When $F$ is a Riemann metric, for any  $C^2$ function $f$ on $M$, the Hessians of $f$ defined by (\ref{hessian1}) and (\ref{hessian2}) respectively are identical.
\end{rem}

\section{Lie derivatives on Finsler manifolds }

Lie derivatives have close connections with Laplacians and Hessians of smooth functions on the manifolds and they are also important tools for studies on Ricci soliton and Ricci flow on Finsler manifolds. However, up to now, there are not yet exact definitions for various Lie derivatives on Finsler manifolds. Further,  some wrong computations about Lie derivatives on Finsler manifolds can be found in some literatures. These cases motivate us to  optimize and improve the various definitions of Lie derivatives on Finsler manifolds.

Let $V=V^{i}\frac{\pa}{\pa x^i}$ be a vector field on $M$ and $\{\varphi _{t}\}$ the local 1-parameter transformation group of $M$ generated by $V$, $V(x)=\frac{d\varphi _{t}(x)}{dt}|_{t=0}$.  The Lie derivative of a tensor in
the direction of $V$ is defined as the first-order term in a suitable Taylor expansion of
the tensor when it is moved by the flow of $V$. The precise formula, however, depends on what type of tensor we use (\cite{Petersen}).

In the following, we mainly consider the Lie derivative on a Finsler manifold $(M, F)$ of dimension $n$. For each $\varphi _{t}$, it is naturally extended to a transformation $\tilde{\varphi}_{t}: \ TM \rightarrow TM$ defined by
\[
\tilde{\varphi}_{t}(x,y):=\left(\varphi _{t}(x), (\varphi _{t})_{*}(y)\right)=\left(\varphi _{t}(x), y^{i}\frac{\pa \varphi _{t}(x)}{\pa x^i}\right).
\]
It is easy to check that $\{\tilde{\varphi}_{t}\}$  is a local 1-parameter transformation group of $TM$. Further,
\beq
\left.\frac{d\tilde{\varphi}_{t}(x, y)}{dt}\right|_{t=0}&=& \left.\left(\frac{d\varphi _{t}(x)}{dt}, y^{m}\frac{\pa ^{2} \varphi _{t}}{\pa t\pa x^{m}}\right)\right|_{t=0}=\left.\left(V(x), y^{m}\frac{\pa}{\pa x^m}\left(\frac{d\varphi _{t}(x)}{dt}\right)\right)\right|_{t=0} \nonumber\\
&=& \left(V(x), y^{m}\frac{\pa V(x)}{\pa x^m}\right).
\eeq
Then $\hat{V}:=V^{i}(x)\frac{\pa}{\pa x^i}+y^{m}\left(\frac{\pa V^{i}}{\pa x^{m}}\right)\frac{\pa}{\pa y^i}$ is the vector field on $TM$ induced by  $\{\tilde{\varphi}_{t}\}$. We call  $\hat{V}$ the complete lift of $V$.

If $f : TM \rightarrow \R$ is a function, then
$f\left(\tilde{\varphi}_{t}(x,y)\right)=f(x,y )+t\left({\cal L}_{\hat{V}} f\right)(x,y)+o(t)$
or
\[
\left({\cal L}_{\hat{V}} f\right)(x,y)=\lim _{t \rightarrow 0} \frac{f\left(\tilde{\varphi}_{t}(x,y)\right)-f(x,y)}{t}.
\]
Thus the Lie derivative ${\cal L}_{\hat{V}}f$ is simply the directional derivative $df(\hat{V})$, that is,
\be
{\cal L}_{\hat{V}} f=\hat{V}f.
\ee

When we have a vetor field $Y\in T(TM)$ things get a litle more complicated as $Y |_{\tilde{\varphi}_{t}(x,y)}$ can't
be compared directly to $Y|_{(x,y)}$ since the vectors live in different tangent spaces. Thus we
consider the curve $t \mapsto \tilde{\varphi}_{t}^{*}\left(Y |_{\tilde{\varphi}_{t}(x, y)}\right)$ that lies in $T_{y} (TM )$, here $\tilde{\varphi}_{t}^{*}:=(\tilde{\varphi}_{t})_{*}^{-1}= (\tilde{\varphi}_{-t})_{*}$.  In other words we define
\[
\left({\cal L}_{\hat{V}} Y\right) |_{(x,y)}=\lim _{t \rightarrow 0} \frac{\tilde{\varphi}_{t}^{*}\left(Y |_{\tilde{\varphi}_{t}(x,y)}\right)-\left.Y\right|_{(x,y)}}{t}.
\]
This Lie derivative turns out to be the Lie bracket (\cite{Petersen}),
\be
{\cal L}_{\hat{V}} Y=[\hat{V}, Y]. \label{LievectorTM}
\ee

If $\xi \in {\cal T}^{0}_{k}(TM)$ is a tensor of $(0, k)-$type over $TM$, its Lie derivative ${\cal L}_{\hat{V}}\xi$ with respect to $\hat{V}$ is the tensor of the same type given by
\be
({\cal L}_{\hat{V}}\xi)(Y_{1}, \cdots, Y_{k}):={\hat{V}}\left(\xi (Y_{1}, \cdots, Y_{k})\right)-\sum\limits_{i=1}^{k}\xi (Y_{1}, \cdots, {\cal L}_{\hat{V}}Y_{i}, \cdots , Y_{k}), \label{0ktensorLie}
\ee
that is,
\be
({\cal L}_{\hat{V}}\xi)(Y_{1}, \cdots, Y_{k})= {\cal L}_{\hat{V}}\left(\xi (Y_{1}, \cdots, Y_{k})\right)-\sum\limits_{i=1}^{k}\xi (Y_{1}, \cdots, {\cal L}_{\hat{V}}Y_{i}, \cdots , Y_{k}),    \label{0ktensorLie}
\ee
where $Y_{i}\in T(TM_{0}), \ 1\leq i \leq k.$

Let $\eta \in {\cal T}^{1}_{k}(TM)$ be a tensor of $(1, k)-$type over $TM$. The Lie derivative ${\cal L}_{\hat{V}}\eta$ of $\eta$ with resect to $\hat{V}$ is defined by (\cite{Lovas})
\be
({\cal L}_{\hat{V}}\eta)(Y_{1}, \cdots, Y_{k}):=\left[{\hat{V}}, \eta (Y_{1}, \cdots, Y_{k})\right]-\sum\limits_{i=1}^{k}\eta (Y_{1}, \cdots, {\cal L}_{\hat{V}}Y_{i}, \cdots , Y_{k}) \label{0ktensorLie}
\ee
for $Y_{i}\in T(TM_{0}), \ 1\leq i \leq k.$ Clearly, ${\cal L}_{\hat{V}}\eta$ is still a tensor of $(1, k)-$type over $TM$.  Obviously, (\ref{0ktensorLie}) can be also rewritten as
\be
({\cal L}_{\hat{V}}\eta)(Y_{1}, \cdots, Y_{k})={\cal L}_{\hat{V}}\left(\eta (Y_{1}, \cdots, Y_{k})\right)-\sum\limits_{i=1}^{k}\eta (Y_{1}, \cdots, {\cal L}_{\hat{V}}Y_{i}, \cdots , Y_{k}).
\ee

\vskip 2mm
As natural applications of the definitions above, when the tensors that we discuss are restricted to the vector bundle $\pi^{*} T M$ and its dual $\pi^{*} T^{*} M$ on $TM$,  we firstly give the following convention:
if $X\in \pi^{*}TM_{0}$, the Lie derivative of $X$ with respect to $\hat{V}$ is given by
\be
{\cal L}_{\hat{V}} X = \rho \left[\hat{V}, X\right]\in \pi^{*}TM_{0}. \label{LovasLie}
\ee
Here, our convention is different from that in \cite{Lovas}.

Based on convention (\ref{LovasLie}), firstly, we give  Lie derivatives of some fundamental tensors on $\pi^{*} T M$ and its dual $\pi^{*} T^{*} M$.
For $y=y^{i}\frac{\pa}{\pa x^i}\in T_{x}M$, let ${\cal Y} = (x, y, y)=y^{i}\pa_{i} |_{(x,y)}\in \pi^{*}TM_{0}$ and  $\xi = {\cal L}({\cal Y})=y_{i}dx^{i}|_{(x,y)}\in \pi^{*}T^{*}M_{0}$, where $y_{i}:=g_{ij}(x,y)y^{j}$. Further,   let
\[
{\cal L}_{\hat{V}}{\cal Y}:=({\cal L}_{\hat{V}}y^{i})\pa _{i}, \ \ {\cal L}_{\hat{V}}\xi :=({\cal L}_{\hat{V}}y_{k})dx^{k}
\]
and
\[
{\cal L}_{\hat{V}}g:=({\cal L}_{\hat{V}}g_{ij})dx^{i}\otimes dx^{j}.
\]
Here,  $g=g_{ij}(x,y)dx^{i}\otimes dx^{j}$ is the the inner product on $\pi^{*}TM_{0}$.
Then, it is easy to show that
\beq
&&{\cal L}_{\hat{V}}y^{i} = 0, \label{Lieyi}\\
&& {\cal L}_{\hat{V}}y_{k} = y^{m}(V_{k|m}+V_{m|k}), \label{Lieyk}\\
&&{\cal L}_{\hat{V}}g_{ij}(x,y) = V_{j|i}(x,y)+V_{i|j}(x,y)+2y^{m}V^{l}_{|m}(x,y)C_{lij}(x,y). \label{gLie}
\eeq
where $V_{i}:=g_{ij}(x,y)V^{j}$. Here we have used the fact that $\frac{\pa V^{m}}{\pa x^{k}}g_{jm}=V_{j|k}-g_{jl}\Gamma ^{l}_{ik}V^{i}$.   Further, the Lie derivative of the Finsler metric $F$ is given by
\be
{\cal L}_{\hat{V}}F={\hat{V}}F = F^{-1}V_{0|0}. \label{LieF}
\ee

More general, the Lie derivative of an arbitrary tensor field $\Upsilon _{I}$ on $\pi^{*}TM$ with respect to the complete lift $\hat{V}$ of a vector field $V=V^{i}(x)\frac{\pa}{\pa x^i}$ on $M$ is defined by
\be
({\cal L}_{\hat{V}}\Upsilon _{I})|_{(x,y)}:=\lim_{t\rightarrow 0}\frac{\tilde{\varphi}^{*}_{t}(\Upsilon _{I}|_{\tilde{\varphi}_{t}(x,y)})-\Upsilon _{I}|_{(x,y)}}{t},
\ee
where $I$ denotes a mixed multi-index. In this case, in the local coordinates $(x^i, y^{i})$ on $TM$, the Lie derivative of an arbitrary mixed tensor field, for example, a tensor $T$ of $(1,2)$-type with the components $T^{i}_{jk}(x,y)$ on $\pi ^{*}TM$, is given by
\be
{\cal L}_{\hat{V}}T:=\left({\cal L}_{\hat{V}}T^{i}_{jk}\right)\pa _{i}\otimes dx^{j}\otimes dx^{k}. \label{Liedetensor}
\ee
Here,
\beqn
{\cal L}_{\hat{V}}T^{i}_{jk}&=&\hat{V}\left(T^{i}_{jk}\right)+\frac{\pa V^{a}}{\pa x^{j}}T^{i}_{ak}+\frac{\pa V^{a}}{\pa x^{k}}T^{i}_{ja}-\frac{\pa V^{i}}{\pa x^{a}}T^{a}_{jk}  \\
&=&V^{m}T^{i}_{jk|m}+y^{m}(V^{l}_{|m})\frac{\pa T^{i}_{jk}}{\pa y^{l}}-T^{m}_{jk}V^{i}_{|m}+T^{i}_{mk}V^{m}_{|j}+T^{i}_{jm}V^{m}_{|k},
\eeqn
where  ``$|$" denotes the horizontal covariant derivative with respect to the Chern connection (see \cite{BM}\cite{JB}).

Now let us define the Lie derivative of Chern connection $D$ with respect to $\hat{V}$.
For $\zeta \in T(TM_{0})$ and $X\in \pi ^{*}TM$, define ${\cal L}_{\hat{V}}D$ by (\cite{Lovas}\cite{Petersen})
\beq
\left({\cal L}_{\hat{V}} D\right)_{\zeta} X &=&\left({\cal L}_{\hat{V}} D\right)(\zeta, X) \nonumber\\
&:=&{\cal L}_{\hat{V}}\left(D_{\zeta} X\right)-D_{{\cal L}_{\hat{V}} \zeta} X - D_{\zeta}({\cal L}_{\hat{V}} X) \nonumber\\
&=& {\cal L}_{\hat{V}}\left(D_{\zeta} X\right)-D_{[\hat{V}, \zeta]} X - D_{\zeta}({\cal L}_{\hat{V}} X). \label{ChernconLie}
\eeq

Write Chern connection 1-form $\omega _{j}^{\ i}=\Gamma ^{i}_{jk}(x,y)dx^{k}$. Put
\be
({\cal L}_{\hat{V}} D)(\pa _{j}, \pa _{k}):= ({\cal L}_{\hat{V}}\Gamma ^{i}_{jk}) \pa _{i}. \label{GammaLie}
\ee

By (\ref{LovasLie}), we have
\be
{\cal L}_{\hat{V}}\pa _{j}= \rho [\hat{V}, \pa _{j}]= -\frac{\pa V^{i}}{\pa x^{j}}\pa _{i}. \label{VpajLIe}
\ee
Further,
\beq
{\cal L}_{\hat{V}}(D_{\pa _{j}}\pa _{k}) &=&{\cal L}_{\hat{V}}(\Gamma ^{i}_{kj}\pa _{i})=({\hat{V}}\Gamma ^{i}_{kj})\pa _{i}-\Gamma ^{m}_{kj}\frac{\pa V^{i}}{\pa x^{m}}\pa _{i}, \label{VDpaij}\\
D_{{\cal L}_{\hat{V}}\pa _{j}}\pa _{k} &=& -\Gamma ^{i}_{kl}\frac{\pa V^{l}}{\pa x^{j}}{\pa _{i}}, \label{VLpajk}\\
D_{\pa _{j}}({\cal L}_{\hat{V}}\pa _{k}) &=& -\left(\frac{\pa ^{2}V^{i}}{\pa x^{k}\pa x^{j}}+\Gamma ^{i}_{mj}\frac{\pa V^{m}}{\pa x^{k}}\right){\pa _{i}}.\label{DpajLV}
\eeq
From (\ref{ChernconLie}) and (\ref{VDpaij})- (\ref{DpajLV}), we obtain
\beqn
({\cal L}_{\hat{V}} D)(\pa _{j}, \pa _{k})&=& {\cal L}_{\hat{V}}(D_{\pa _{j}}\pa _{k})- D_{{\cal L}_{\hat{V}}\pa _{j}}\pa _{k} - D_{\pa _{j}}({\cal L}_{\hat{V}}\pa _{k})\\
&=& \left(({\hat{V}}\Gamma ^{i}_{kj})-\Gamma ^{m}_{kj}\frac{\pa V^{i}}{\pa x^{m}}+ \Gamma ^{i}_{kl}\frac{\pa V^{l}}{\pa x^{j}}+\frac{\pa ^{2}V^{i}}{\pa x^{k}\pa x^{j}}+\Gamma ^{i}_{mj}\frac{\pa V^{m}}{\pa x^{k}}\right)\pa _{i}.
\eeqn

From (\ref{GammaLie}), we obtain
\be
{\cal L}_{\hat{V}} \Gamma_{j k}^{i}=\frac{\partial^{2} V^{i}}{\partial x^{j} \partial x^{k}}+\Gamma_{lj}^{i} \frac{\partial V^{l}}{\partial x^{k}}+\Gamma_{kl}^{i} \frac{\partial V^{l}}{\partial x^{j}}-\Gamma_{kj}^{l} \frac{\partial V^{i}}{\partial x^{l}}+V^{l} \frac{\partial \Gamma_{kj}^{i}}{\partial x^{l}}+y^{s} \frac{\partial V^{l}}{\partial x^{s}} P_{k \ jl}^{\ i}, \label{LieGamma1}
\ee
where $P_{k \ jl}^{\ i}:=\frac{\pa \Gamma ^{i}_{kj}}{\pa y^{l}}$ determine the Landsberg curvature of $F$.

Note that
\be
V^{i}_{|j}=\frac{\delta V^{i}}{\delta x^j}+V^{m}\Gamma ^{i}_{mj}=\frac{\pa V^i}{\pa x^j}+V^{m}\Gamma ^{i}_{mj}. \label{Vik}
\ee
We get
\beqn
\frac{\delta V^{i}_{|j}}{\delta x^{k}}&=& \frac{\delta}{\delta x^k}\left(\frac{\pa V^i}{\pa x^j}\right)+ \frac{\delta V^m}{\delta x^k}\Gamma ^{i}_{mj}+V^{m}\frac{\delta \Gamma ^{i}_{mj}}{\delta x^k}\\
&=& \frac{\pa ^{2}V^{i}}{\pa x^{k}\pa x^{j}}+ \frac{\pa V^m}{\pa x^k}\Gamma ^{i}_{mj}+V^{m}\frac{\delta \Gamma ^{i}_{mj}}{\delta x^k}.
\eeqn
Further, we have
\beq
V^{i}_{|j|k}&=& \frac{\delta V^{i}_{|j}}{\delta x^k}+V^{m}_{|j}\Gamma ^{i}_{mk}-V^{i}_{|m}\Gamma ^{m}_{jk} \nonumber\\
&=&\frac{\pa ^{2}V^{i}}{\pa x^{k}\pa x^{j}}+V^{m}\frac{\delta \Gamma ^{i}_{mj}}{\delta x^k}+V^{a}\Gamma ^{m}_{aj}\Gamma ^{i}_{mk}-V^{a}\Gamma ^{i}_{am}\Gamma ^{m}_{jk}\nonumber\\
&& +\frac{\pa V^m}{\pa x^k}\Gamma ^{i}_{mj}+\frac{\pa V^m}{\pa x^j}\Gamma ^{i}_{mk}-\frac{\pa V^i}{\pa x^m}\Gamma ^{m}_{jk}. \label{Vikj}
\eeq
Then, by comparing (\ref{LieGamma1}) and (\ref{Vikj}), we can get
\[
V_{|j| k}^{i}=\mathcal{L}_{\hat{V}} \Gamma_{j k}^{i}-\frac{\partial V^{l}}{\partial x^{s}} y^{s} P_{j \ k l}^{\ i}+ V^{m} \frac{\delta \Gamma_{m j}^{i}}{\delta x^{k}}+V^{r} \Gamma_{r j}^{l} \Gamma_{l k}^{i}-V^{r} \Gamma_{r l}^{i} \Gamma_{j k}^{l}-V^{l}\frac{\pa \Gamma ^{i}_{kj}}{\pa x^l}.
\]
On the other hand,
\beqn
 R_{j \ k m}^{\ i} V^{m} &=& V^{m}\left(\frac{\delta \Gamma_{m j}^{i}}{\delta x^{k}}-\frac{\delta \Gamma_{j k}^{i}}{\delta x^{m}}+\Gamma_{j m}^{l} \Gamma_{k l}^{i}-\Gamma_{l m}^{i} \Gamma_{j k}^{l}\right) \\
 &=& V^{m}\left(\frac{\partial \Gamma_{m j}^{i}}{\partial x^k}-N_{k}^{r} \frac{\partial \Gamma_{m j}^{i}}{\partial y^{r}}-\frac{\partial \Gamma_{j k}^{i}}{\partial x^m}+N_{m}^{r} \frac{\partial \Gamma_{j k}^{i}}{\partial y^{r}}+\Gamma_{j m}^{l} \Gamma_{k l}^{i}-\Gamma_{l m}^{i} \Gamma_{j k}^{l}\right).
\eeqn
Here, $R^{\ i}_{j \ km}=R^{\ i}_{j \ km}(x,y)$ are the coefficients of Riemann curvature tensor with respect to Chern connection. Thus we can get
\beq
 V_{|j| k}^{i}- R_{j \ k m}^{\ i} V^{m}&=& {\cal L}_{\hat{V}} \Gamma_{j k}^{i}-\frac{\partial V^{l}}{\partial x^{m}} y^{m} P_{j \ k l}^{\ i}-V^{m} N_{m}^{r} \frac{\partial \Gamma_{j k}^{i}}{\partial y^{r}}\nonumber \\
 &=& {\cal L}_{\hat{V}} \Gamma_{j k}^{i}-\left(\frac{\partial V^{r}}{\partial x^{s}} +V^{m} \Gamma_{m s}^{r} \right)y^{s} P_{j \ kr}^{\ i}  \nonumber\\
 &=&{\cal L}_{\hat{V}} \Gamma_{j k}^{i}-y^{m} V_{ | m}^{r} P_{j \ k r}^{\ i}.
\eeq

Hence, we have the following
\begin{thm}
The the Lie derivative of  Chern connection with respect to the complete lift $\hat{V}$ of a vector field $V=V^{i}(x)\frac{\pa}{\pa x^i}$ on $M$ is determined by
\be
{\cal L}_{\hat{V}}\Gamma^{i}_{jk}=R^{\ i}_{j \ mk}V^{m}+V^{i}_{|j|k}+y^{m}V^{l}_{|m} P^{\ i}_{j \ kl}. \label{GammaLieD}
\ee
\end{thm}

\begin{rem}
Let $G = y^{i}\frac{\pa}{\pa x^{i}}-2 G^{i}\frac{\pa}{\pa y^{i}}$ be the spray induced by Finsler metric $F$, where $G^{i}(x,y)=\frac{1}{2}\Gamma ^{i}_{jk}(x,y)y^{j}y^{k}$.  By (\ref{LievectorTM}), we have
\be
{\cal L}_{\hat{V}}G  = [\hat{V}, G] = -y^{i}y^{j}\left\{\frac{\pa ^{2}V^{k}}{\pa x^{i}\pa x^{j}}+V^{m}\frac{\pa \Gamma ^{k}_{ij}}{\pa x^{m}}+2\frac{\pa V^{m}}{\pa x^{j}}\Gamma ^{k}_{mi}-\Gamma ^{m}_{ij}\frac{\pa V^{k}}{\pa x^{m}}\right\}\frac{\pa}{\pa y^{k}}.\label{Liespary}
\ee
Let
\[
{\cal L}_{\hat{V}}G:= ({\cal L}_{\hat{V}}G^{k})\frac{\pa}{\pa y^{k}}.
\]
Comparing (\ref{LieGamma1}) with (\ref{Liespary}), we have
\beq
{\cal L}_{\hat{V}}G^{k}&=& -y^{i}y^{j}\left\{\frac{\pa ^{2}V^{k}}{\pa x^{i}\pa x^{j}}+V^{m}\frac{\pa \Gamma ^{k}_{ij}}{\pa x^{m}}+2\frac{\pa V^{m}}{\pa x^{j}}\Gamma ^{k}_{mi}-\Gamma ^{m}_{ij}\frac{\pa V^{k}}{\pa x^{m}}\right\}\nonumber\\
&=& -y^{i}y^{j}{\cal L}_{\hat{V}}\Gamma^{k}_{ij}. \label{LieGGamma}
\eeq
Further, by (\ref{GammaLieD}), we know the following
\be
{\cal L}_{\hat{V}}G^{k}= -(R^{k}_{\ m}V^{m} + V^{k}_{|0|0}).
\ee
Obviously, ${\cal L}_{\hat{V}}G^{k}\neq \frac{1}{2}\left({\cal L}_{\hat{V}}\Gamma^{k}_{ij}\right)y^{i}y^{j}$.
\end{rem}

\vskip 3mm
It  is natural to establish a connection between Lie derivative and Hessian (or Laplacian) of a $C^{2}$ function on the manifold. Firstly,
given a $C^2$ function $u=u(x)$ on $M$, let us determine the Lie derivative of the fundamental tensor with respect to the complete lift $\widehat{\nabla u}$ of the gradient vector field $\nabla u$. Recall that $\nabla u=g^{ij}(x,\nabla u)\frac{\pa u}{\pa x^i}\frac{\pa}{\pa x^j}$ and
\beqn
&& \nabla ^{i}u(x)=g^{ij}(x,\nabla u)\frac{\pa u}{\pa x^j}= g^{*ij}(x, du)\frac{\pa u}{\pa x^j},\\
&& (\nabla u)_{i}=g_{ij}(x, \nabla u)(\nabla ^{j}u)= u_{|i}(x).
\eeqn
The complete lift $\widehat{\nabla u}$ of $\nabla u$ is given by $\widehat{\nabla u}=(\nabla ^{i}u)(x)\frac{\pa}{\pa x^i}+y^{m}\left(\frac{\pa (\nabla^{i}u)}{\pa x^{m}}\right)\frac{\pa}{\pa y^i}$.
Let
\[
(\widehat{\nabla u})_{i}:=g_{im}(x,y) (\nabla ^{m}u)(x)=g_{im}(x,y)\left(g^{mk}(x,\nabla u)\frac{\pa u}{\pa x^k}\right).
\]
By (\ref{gLie}), we have
\be
{\cal L}_{\widehat{\nabla u}}g_{ij}(x,y)=(\widehat{\nabla u})_{i|j}(x,y)+(\widehat{\nabla u})_{j|i}(x,y)+2y^{m}(\nabla ^{l}u)_{|m}(x,y)C_{lij}(x,y). \label{GiLiexy}
\ee
Here, by the definitions, we have
\beqn
(\widehat{\nabla u})_{i|j}(x,y)&=& g_{im}(x,y)\left(g^{mk}(x,\nabla u)\right)_{|j}(x,y)\frac{\pa u}{\pa x^k} \\
&& +g_{im}(x,y)g^{mk}(x,\nabla u)u_{|k|j}(x,y) \\
(\nabla ^{l}u)_{|m}(x,y)&=& \left(g^{lr}(x,\nabla u)\frac{\pa u}{\pa x^r}\right)_{|m}(x,y)=\left(g^{lr}(x,\nabla u)\right)_{|m}(x,y)\frac{\pa u}{\pa x^r}\\
&& +g^{lr}(x,\nabla u)u_{|r|m}(x,y).
\eeqn
Further,
\beqn
\left(g^{mk}(x,\nabla u)\right)_{|j}(x,y)&=& -g^{mr}(x,\nabla u)g^{ks}(x,\nabla u)\left(g_{rs}(x,\nabla u)\right)_{|j}(x,y)\\
&=&-g^{mr}(x,\nabla u)g^{ks}(x,\nabla u)\Big\{\frac{\pa g_{rs}}{\pa x^j}(x,\nabla u)+2C_{rsl}(x,\nabla u)\frac{\pa (\nabla ^{l}u)}{\pa x^{j}} \\
&&-g_{ls}(x,\nabla u)\Gamma ^{l}_{rj}(x,y)-g_{rl}(x,\nabla u)\Gamma ^{l}_{sj}(x,y)\Big\} \\
&=&-g^{mr}(x,\nabla u)g^{ks}(x,\nabla u)\Big\{\frac{\pa g_{rs}}{\pa x^j}(x,\nabla u)+2C_{rsl}(x,\nabla u)\frac{\pa (\nabla ^{l}u)}{\pa x^{j}}\Big\}\\
&&+g^{mr}(x,\nabla u)\Gamma ^{k}_{rj}(x,y)+g^{ks}(x,\nabla u)\Gamma ^{m}_{sj}(x,y).
\eeqn
In particular, it is easy to see that
\[
(\widehat{\nabla u})_{i|j}(x,\nabla u)=u_{|i|j}(x,\nabla u).
\]
Thus, by (\ref{GiLiexy}), we get
\be
{\cal L}_{\widehat{\nabla u}}g_{ij}(x,\nabla u)=2u_{|i|j}(x,\nabla u)+2(\nabla ^{m}u)u_{|r|m}(x,\nabla u)C^{r}_{ij}(x, \nabla u). \label{gijGVL}
\ee
Then, by (\ref{hessian}), we know that the Hessian in Finsler geometry can also be determined by Lie derivative.
\begin{prop}
Let $u=u(x)$ be a  $C^2$ function on Finsler manifold $(M, F)$. Then the Hessian of $u$ is determined by
\be
{\rm Hess}(u)=\frac{1}{2}{\cal L}_{\widehat{\nabla u}}g(x,\nabla u)-\left((\nabla ^{m}u)u_{|r|m}(x,\nabla u)C^{r}_{ij}(x, \nabla u)\right),
\ee
where $g(x,\nabla u):=\left(g_{ij}(x,\nabla u)\right)$.
\end{prop}
\begin{cor}{\rm (\cite{JW})}
Let $u=u(x)$ be a  $C^2$ function on Riemannian manifold $(M, g)$. Then the Hessian of $u$ is determined by
\[
{\rm Hess}(u)=\frac{1}{2}{\cal L}_{\nabla u}g.
\]
\end{cor}

\section{Gradient estimate on Randers manifolds}

Let $F=\alpha +\beta$ be a Randers metric on a manifold $M$ with $\|\beta\|_{\alpha}<1$. Take a basis $\{{\bf b}_{i}\}^{n}_{i=1}$ for $TM$ and its dual basis $\{\theta ^{i}\}_{i=1}^{n}$ for $T^{*}M$.  Express
$\alpha$ and $\beta$ by
\[
\alpha (x,y)=\sqrt{a_{ij}(x)y^{i}y^{j}}, \ \ \ \beta (x,y)=b_{i}(x)y^{i}, \ \ y=y^{i}{\bf b}_{i}.
\]
The dual metric $F^{*}$ of $F$ is still of Randers type  on $T^{*}M$. More precisely, we can get $F^{*}(x, \xi)=\alpha ^{*}(x, \xi)+\beta ^{*}(x, \xi)$, where
\[
\alpha ^{*}(x, \xi)=\sqrt{a^{*ij}(x)\xi _{i}\xi _{j}}, \ \ \ \beta ^{*}(x, \xi) =b^{*i}(x)\xi _{i}, \ \ \xi =\xi _{i} \theta ^{i}
\]
and
\beq
a^{*ij}&=&\frac{(1-\|\beta\|_{\alpha}^{2})a^{ij}+b^{i}b^{j}}{(1-\|\beta\|_{\alpha}^{2})^2},\\
b^{*i}&=&-\frac{b^{i}}{1-\|\beta\|_{\alpha}^{2}},
\eeq
where $b^{i}:=b_{j}a^{ij}$ (see \cite{HrSh}). Equivalently, we have
\beq
a^{ij}&=&(1-\|\beta\|_{\alpha}^{2})(a^{*ij}-b^{*i}b^{*j}),\\
b^{i}&=&-(1-\|\beta\|_{\alpha}^{2})b^{*i}.
\eeq
Let $(a_{*ij}):=(a^{*ij})^{-1}$ and $b_{*i}:=a_{*ij}b^{*j}$. Then
\beq
a_{*ij}&=&(1-\|\beta\|_{\alpha}^{2})(a_{ij}-b_{i}b_{j}),\\
b_{*i}&=&-(1-\|\beta\|_{\alpha}^{2})b_{i}.
\eeq
Equivalently,
\beq
a_{ij}&=&\frac{(1-\|\beta\|_{\alpha}^{2})a_{*ij}+b_{*i}b_{*j}}{(1-\|\beta\|_{\alpha}^{2})^2}, \label{aijastar}\\
b_{i}&=&-\frac{b_{*i}}{1-\|\beta\|_{\alpha}^{2}}.\label{bistar}
\eeq
The norm $\|\beta ^{*}\|_{\alpha ^{*}}:=\sup_{\alpha ^{*}(\xi)=1}\beta ^{*}(\xi)$ is given by
\[
\|\beta ^{*}\|_{\alpha ^{*}}=\sqrt{a_{*ij}b^{*i}b^{*j}}=\sqrt{\frac{1}{1-\|\beta\|_{\alpha}^{2}}(a_{ij}-b_{i}b_{j})b^{i}b^{j}}=\|\beta\|_{\alpha}.
\]
Then, by putting $b^{*}:=\|\beta ^{*}\|_{\alpha ^{*}}$ and $b:=\|\beta\|_{\alpha}$, we have $b^{*}=b$.

Now let us consider the Legendre transformations ${\cal L}: TM \rightarrow T^{*}M$  on Finsler manifold $(M, F=\alpha +\beta)$. Let $\xi ={\cal L}(y)=\xi _{i}\theta ^{i}$. Then
\be
\xi _{i}=g_{ij}(x,y)y^{j}=F(x,y)\left(\frac{a_{ij}y^{j}}{\alpha (x,y)}+b_{i}\right). \label{xiiR}
\ee
Conversely, Let $y={\cal L}^{-1}(\xi)=y^{i}{\bf b}_{i}$. We have
\be
y^{i}=g^{*il}(x,\xi)\xi _{l}=F^{*}(x,\xi)\left(\frac{a^{*il}\xi _{l}}{\alpha ^{*}(x,\xi)}+b^{*i}\right). \label{yiR}
\ee
By (\ref{aijastar}), (\ref{bistar}) and (\ref{yiR}), we obtain
\[
\alpha ^{2}(x,y)=a_{ij}(x)y^{i}y^{j}=\left(\frac{F^{*2}(x,\xi)}{\alpha ^{*}(x,\xi)}\right)^{2}\frac{1}{(1-b^2)^2},
\]
that is,
\be
\alpha (x,y)=\left(\frac{1}{1-b^2}\right)\frac{F^{*2}(x,\xi)}{\alpha ^{*}(x,\xi)}  \label{alphastar}
\ee
and
\be
\beta (x,y)=b_{i}y^{i}=-\frac{F^{*}(x,\xi)}{1-b^2}\left(\frac{\beta ^{*}(x,\xi)}{\alpha ^{*}(x,y)}+b^{2}\right).\label{betastar}
\ee
Then, by using (1.57) in \cite{ChernShen} and by (\ref{alphastar}), (\ref{betastar}), we obtain
\beq
g^{*ij}(x, \xi)&=&g^{ij}(x,y)=\frac{\alpha}{F}a^{ij}+\left(\frac{\alpha}{F}\right)^{2}\frac{\beta +\alpha b^{2}}{F}\frac{y^i}{\alpha}\frac{y^j}{\alpha}-\left(\frac{\alpha}{F}\right)^{2}\left(b^{j}\frac{y^i}{\alpha}+b^{i}\frac{y^j}{\alpha}\right)\nonumber\\
&=&\frac{F^{*}}{\alpha ^{*}}(a^{*ij}+b^{*i}b^{*j})+\frac{F^{*}}{\alpha ^{*2}}(\xi ^{i}b^{*j}+\xi ^{j}b^{*i}) \nonumber \\
&&-\frac{\beta ^{*}}{\alpha ^{*}}\left(\frac{\xi ^{i}}{\alpha ^{*}}+b^{*i}\right)\left(\frac{\xi ^{j}}{\alpha ^{*}}+b^{*j}\right), \label{gijstar}
\eeq
where $\xi ^{i}:=a^{*ij}\xi _{j}$.

In the following, we will determine the gradient vector field $\nabla f$ of  a smooth function $f$ on a Randers manifold $(M, F=\alpha +\beta)$. By the definition,  $\nabla f(x)=\nabla ^{i}f\frac{\pa}{\pa x^{i}}$, where $\nabla ^{i}f:= g^{ij}(x,\nabla f)\frac{\pa f}{\pa x^j}= g^{*ij}(x, df)\frac{\pa f}{\pa x^j}$ on $M_{f}$. Put
\[
\beta ^{\sharp}:=b^{i}(x)\frac{\pa}{\pa x^{i}}, \ \ \nabla ^{\alpha}f:=(df)^{\sharp}=a^{ij}(x)\frac{\pa f}{\pa x^j}\frac{\pa}{\pa x^{i}}.
\]
Let $\langle \beta ^{\sharp},\nabla ^{\alpha}f\rangle =b^{i}(x)\frac{\pa f}{\pa x^i}$ denote the inner product of $ \beta ^{\sharp}$ and $\nabla ^{\alpha}f$ with respect to $\alpha$. Further, put
\[
A:=\sqrt{(1-b^2)\|\nabla ^{\alpha}f\|^{2}_{\alpha}+\langle \beta ^{\sharp},\nabla ^{\alpha}f\rangle^{2}},
\]
where $\|\nabla ^{\alpha}f\|_{\alpha}$ denotes the norm of the gradient vector field $\nabla ^{\alpha}f$ with respect to ${\alpha}$. Then we can get the following
\beq
\alpha ^{*}(x,df)&=&\sqrt{a^{*ij}(x)\frac{\pa f}{\pa x^i}\frac{\pa f}{\pa x^j}}=\frac{A}{1-b^2}, \label{alphadf}\\
\beta ^{*}(x,df)&=& -\frac{\langle \beta ^{\sharp},\nabla ^{\alpha}f\rangle}{1-b^2}, \label{betadf} \\
F^{*}(x,df)&=&\frac{A-\langle \beta ^{\sharp},\nabla ^{\alpha}f\rangle}{1-b^2}, \label{Fdf}\\
\xi ^{j}\frac{\pa f}{\pa x^j}&=&a^{*jl}(x)\frac{\pa f}{\pa x^l}\frac{\pa f}{\pa x^j}=\alpha ^{*2}(x,df)=\left(\frac{A}{1-b^2}\right)^{2}. \label{xidf}
\eeq
Here, $\xi ^{j}=a^{*jl}\frac{\pa f}{\pa x^l}$. Thus, by (\ref{gijstar}), we can obtain
\be
\nabla ^{i}f= g^{*ij}(x, df)\frac{\pa f}{\pa x^j}=\frac{A-\langle \beta ^{\sharp},\nabla ^{\alpha}f\rangle}{A(1-b^2)}a^{ij}(x)\frac{\pa f}{\pa x^j}-\frac{(A-\langle \beta ^{\sharp},\nabla ^{\alpha}f\rangle)^2}{A(1-b^2)^2}b^{i}.\label{nablaf}
\ee
\begin{thm}\label{Randersgrad}
The gradient vector field $\nabla f$ of a smooth function $f$ on Randers manifold $(M, F=\alpha +\beta)$ is given by
\be
\nabla f =\frac{A-\langle \beta ^{\sharp},\nabla ^{\alpha}f\rangle}{A(1-b^2)}\nabla ^{\alpha}f-\frac{(A-\langle \beta ^{\sharp},\nabla ^{\alpha}f\rangle)^2}{A(1-b^2)^2}\beta ^{\sharp}. \label{gradientRanders}
\ee
Further, we have the following gradient estimate formula on Randers manifold  $(M, F=\alpha +\beta)$
\be
F(x,\nabla f)\leq \frac{\alpha (x,\nabla ^{\alpha}f)}{1-b}. \label{gradientestimate}
\ee
\end{thm}
{\bf Proof.}   It is easy to see that (\ref{gradientRanders}) holds by (\ref{nablaf}). By Cauchy-Schwarz inequality, we know that
\[
A=\sqrt{(1-b^2)\|\nabla ^{\alpha}f\|^{2}_{\alpha}+\langle \beta ^{\sharp},\nabla ^{\alpha}f\rangle^{2}}\leq \sqrt{(1-b^2)\|\nabla ^{\alpha}f\|^{2}_{\alpha}+b^{2}\|\nabla ^{\alpha}f\|^{2}_{\alpha}}=\|\nabla ^{\alpha}f\|_{\alpha}.
\]
Then
\[
A-\langle \beta ^{\sharp},\nabla ^{\alpha}f\rangle \leq \|\nabla ^{\alpha}f\|_{\alpha}+b \|\nabla ^{\alpha}f\|_{\alpha}=(1+b)\|\nabla ^{\alpha}f\|_{\alpha}.
\]
Hence, by (\ref{Fdf}), we have
\beqn
F(x,\nabla f)&=&F^{*}(x,df)=\frac{A-\langle \beta ^{\sharp},\nabla ^{\alpha}f\rangle}{1-b^2}\\
&\leq & \frac{\|\nabla ^{\alpha}f\|_{\alpha}}{1-b}=\frac{\alpha (x,\nabla ^{\alpha}f)}{1-b}.
\eeqn
\qed

Theorem \ref{Randersgrad} shows that the gradient vector field $\nabla f$ of a smooth function $f$ on Randers manifold $(M, F=\alpha +\beta)$ is determined completely by the gradient vector field $\nabla ^{\alpha}f$ of  $f$ with respect to Riemann metric $\alpha$.

\vskip 8mm

\noindent
Xinyue Cheng \\
School of Mathematical Sciences \\
Chongqing Normal University \\
Chongqing  401331,  P. R. of China  \\
E-mail:  chengxy@cqnu.edu.cn

\end{document}